\documentclass[11pt,oneside,a4paper]{article}
\usepackage{graphicx}
\usepackage{amsmath}
\usepackage{enumerate}
\begin{document}

\title{New exact solutions of nonlinear variants of the RLW, the PHI-four
and Boussinesq equations based on modified extended direct algebraic method}

\author{A. A. Soliman\footnote{Suez Canal University, Faculty of Education (AL-Arish), Department of Mathematics, AL-Arish 45111, Egypt, asoliman\_99@yahoo.com }\qquad \qquad  H. A. Abdo\footnote{Freie Universit\"{a}t Berlin, Institut f\"{u}r Informatik, Takustra\ss e 9, D-14195 Berlin, Germany,  abdo@inf.fu-berlin.de } }

\maketitle
\begin{abstract}
By means of modified extended direct algebraic method (MEDA) the multiple
exact complex solutions of some different kinds of nonlinear partial
differential equations are presented and implemented in a computer algebraic
system. New complex solutions for nonlinear equations such as the variant of
the RLW equation, the variant of the PHI-four equation and the variant
Boussinesq equations are obtained
\end{abstract}

\section{Introduction}

\noindent \qquad It is well known that the nonlinear physical phenomena are related to nonlinear partial
 differential equations which are involved in many fields such as physics, chemistry, mechanics, etc.
 As mathematical models of the phenomena, the investigation of exact solutions of the partial
 differential equations will help one to understand these phenomena better. In recent years, many
 powerful methods to construct exact solutions of nonlinear partial differntial equations have been
 established and developed. Among these are variational iteration method
 \cite{abdou:vim,abdou:variational, darvishi:numerical,soliman:numerical,soliman:2dburger, soliman:nuvim},
 tanh function method \cite{fan:soliton,parker:automated}, modified extended tanh function method
 \cite{abdou:modi,ali:newexact,wakil:modified,soliman:met,soliman:exacttrav}, sine-cosine method
 \cite{ali:solution,yan:trans} Exp-method\cite{abdou:newapp}, inverse scattering method \cite{gardner:kdv},
 Hirota's bilinear method \cite{hirota:exact}the homogeneous balance method\cite{wang:exactkdv},
the Riccati expansion method with constant coefficients \cite{yan:explicit,yan:new}. The regularized long wave
(RLW) equation is an important nonlinear wave equation. Solitary waves are wave packet or pulses, which
propagate in nonlinear dispersive media. Due to dynamical balance between the non-linear and dispersive effects
these waves retain a stable wave form. The regularized long wave (RLW) equation is an alternative description of
non-linear dispersive waves to the more usual Kortewege-de vries (KDV) equation \cite{peregrine:calc}. Numerical
solutions based on finite difference techniques \cite{abdullove:soliton,eilbek:numerical}. Rung-Kutta method
\cite{bona:numerical} and Galerkin's method ~\cite{alexander:galerkin} have been given. Alexander and Morris
\cite{alexander:galerkin} constructed a global trial function mainly from cubic splines. Gardner and Gardner
\cite{gardner:solitary}, using the Galerkin's method and cubic B-spline as element shape function to construct
an implicit finite element solution. The least squares method using linear space-time finite elements used to
solve the RLW equation \cite{gardner:least}, Soliman and Raslan \cite{soliman:collocation} solved the RLW
equation by using collocation method using quadratic B-spline at the mid point. Soliman and Hussien
\cite{soliman:rlw} used the collocation method with septic spline to solve the RLW equation, Soliman
\cite{soliman:numerical} using the finite difference method with the similarity solution of the partial
differential equations to obtain the numerical scheme for RLW equation this approach eliminate the difficult
associated to the boundary. in the end Soliman \cite{soliman:numerical} using the variational iteration method
to obtain the exact solutions of the gneralized RLW. The PHI-four equation is considered as a particular form of
the Klein-Gordon equation that model phenomenon in particle physics where kink and anti-kink solitary waves
interact\cite{benjamin:model}. Recently, the direct algebraic method and symbolic computation have been
suggested to obtain the exact complex solutions of nonlinear partial differential equations
\cite{soliman:medam,zhang:dam,zhang:netws}.The aim of this paper is to extend the modified extended direct
algebraic method MEDA method to solve three different types of nonlinear differential equations such as the
variant of the RLW, the variant ~of the PHI-four and variant Boussinesq equations
\cite{ali:newexact,soliman:exacttrav}.
\section{\protect\bigskip \textbf{Modified extended direct algebraic method}}

To illustrate the basic concepts of the modified extended direct algebraic
(MEDA) method. We consider a given PDE in two independent variables given
by
\begin{equation}
\label{eq:m201}
       F(u,u_{x},u_{t},u_{xx},....)= 0
\end{equation}

We first consider its travelling solutions $u(x,t) = u(z)$, $z=i(x + c\:t)$ \
or $z=i(x - c\:t)$, $i=\sqrt{-1}$, then Eq.~\eqref{eq:m201} becomes an ordinary differential
equation
\bigskip 
\begin{equation}
\label{eq:m202}
   H(u,i\:u^{\prime },-\:i\:c\:u^{\prime },-\:u^{\prime \prime },...)=0
\end{equation}
where $u^{\prime }=\frac{du}{dz}.$
In order to seek the solutions of Eq.~\eqref{eq:m201}, we introduce the following ansatze
\begin{equation}
\label{eq:m203}
  u(z)=a_{0}+\sum_{j=1}^{M}(a_{j}\:\phi ^{j}+ b_{j}\:\phi ^{-j})
\end{equation}
\begin{equation}
\label{eq:m204}
       \phi ^{\prime }= b + \phi ^{2}
\end{equation}

where $b$ is a parameter to be determined, $\phi =\phi(z)$, $\phi ^{\prime
}=\frac{d\phi}{dz}$. \ The parameter $M$ can be found by balancing the
highest-order derivative term with the nonlinear terms \cite{wang:solitary}.
Substituting~\eqref{eq:m203} into~\eqref{eq:m202} with~\eqref{eq:m204} will
yield a system of algebraic equations with respect to $a_{j}$, $b_{j}$, $b$ and $c$(where $j=1...M$) because all the coefficients of $\phi ^{j}$ have to vanish. we can then
determine $a_{0}$, $a_{j}$, $b_{j}$, $b$, and $c$. Eq.~\eqref{eq:m204} has the general
solutions:

\begin{enumerate}[(I)]
	\item \textbf{If $b<0$}
	\begin{equation*}
      \phi =-\sqrt{-b}\:\tanh\:(\sqrt{-b}\:z),\ \ \ or\ \ \  \phi =-\sqrt{-b}\:\coth\:(\sqrt{-b}\:z)
  \end{equation*}
it depends on initial conditions.
	\item \textbf{If $b>0$}
  \begin{equation*}
      \phi =\sqrt{b}\:\tan\:(\sqrt{b}\:z),\ \ \ or\ \ \  \phi =-\sqrt{b}\:\cot\:(\sqrt{b}\:z)
  \end{equation*}
it depends on initial conditions.
  \item \textbf{If $b=0$}
  \begin{equation}
  \label{eq:m205}
      \phi =\frac{-1}{z}
  \end{equation}
\end{enumerate}
Substituting the results into~\eqref{eq:m203}, then we obtain the exact travelling wave
solutions of Eq.~\eqref{eq:m201}.
To illustrate the procedure, three examples related to variant of the RLW,
variant of the PHI-four and variant Boussinesq equations are given in the
following.
\section{Applications}
\subsection{\protect\bigskip Variant of the regularized long-wave (RLW) equation}

\bigskip Let us first consider the nonlinear variant of the regularized long
wave equation which has the form~\cite{soliman:exacttrav}
\begin{equation}
\label{eq:m3st1}
      u_{t}+\alpha \:u_{x}-\lambda \:(u^{n})_{x}+\beta \:(u^{n})_{xxt} = 0
\end{equation}

where $\alpha$, $\lambda$, and $\beta$ are arbitrary constants. In order
to solve Eq.~\eqref{eq:m3st1} by the MEDA method, we use the wave transformation
$u(x,t)=U(z)$ with wave complex variable $z=i(x-c\:t),$ Eq.~\eqref{eq:m3st1} takes the form
of an ordinary differential equation as

\begin{equation}
\label{eq:m3st2}
    (\alpha -c)\:U^{\prime }-\lambda\:(U^{n})^{\prime }+\beta \:c\:(U^{n})^{\prime
     \prime \prime}=0
\end{equation}

\bigskip Integrating Eq.~\eqref{eq:m3st2} once with respect to $z$ and setting the
constant of integration to be zero, we obtain

\begin{equation}
\label{eq:m3st3}
      (\alpha - c)\:U -\lambda \:U^{n}+\beta \: c \:(U^{n})^{\prime \prime }= 0
\end{equation}

Or equivalently

\begin{equation}
\label{eq:m3st4}
    (\alpha -c)\:U-\lambda \: U^{n}+\beta \: c\:n\:U^{n-1}U^{\prime \prime }+\beta
     \:c\:n\:(n-1)\:U^{n-2}\:(U^{\prime })^{2}=0
\end{equation}

Balancing the order of $U^{n}$ with the order of $U^{n-1}\:U^{\prime \prime }$ in
Eq.~\eqref{eq:m3st4}, we find $M = -\frac{2}{n-1}$. To get a closed form analytic solution, 
the parmeter $M$ should be an integer. A transformation formula 
$U = V^{-\frac{1}{n-1}}$ should be used to obtain this analytic solution. So 
Eq.~\eqref{eq:m3st4} takes the form

\begin{equation}
\label{eq:m3st5}
(\alpha - c)(n-1)^{2}\:V^{3}-\lambda\:(n-1)^{2}\:V^{2}-\beta c\:n\:(n-1)V\:V^{\prime
\prime}+\beta\:c\:n(2n-1)\:(V^{\prime })^{2}=0
\end{equation}

Balancing the order of $V^{3}$ with the order of $V\:V^{\prime \prime }$ in
Eq.~\eqref{eq:m3st5}, gives $M = 2$. So the solution takes the form

\begin{equation}
\label{eq:m3st6}
    V(z)=a_{0}+a_{1}\: \phi(z)+ a_{2}\: \phi(z)^{2}+ b_{1}\: \phi(z)^{-1}+ 
    b_{2}\: \phi \:(z)^{-2}
\end{equation}

Inserting Eq.~\eqref{eq:m3st6} into Eq.~\eqref{eq:m3st5} and making use of 
Eq.~\eqref{eq:m204}, using the Maple Package, we get a system of algebraic equations, for
$a_{0},\:a_{1},\:a_{2},\:b_{1},\:b_{2}$ and $b$. We solve the obtained system of
algebraic equations give the following three cases:

\bigskip Case (I): consider, $A = (-n \: \lambda -\lambda + 2 \: \alpha \: n \: a_{0})$ then,
         $a_{1} = 0,\; \; a_{2} = 0,\; \; b_{1} = 0 $,\; \;
         $b_{2}=\frac{a_{0}^{2}\: \lambda \:(n^{2}-2 \: n+1)}{2 \: n \: \beta \: A}$,\; \;  
         $c=\frac{A}{2 \: n \: a_{0}}$,
         \; \; $b = \frac{a_{0} \: \lambda \:(n^{2}-2 \: n+1)}{2 \: n \: \beta \: A}$ , with $a_{0}$ being an arbitrary constant.

\begin{equation}
\label{eq:m3st7}
       V(x,t)= a_{0}( 1  +  \cot^{2}  (\frac{i}{2}  \sqrt{\frac{
2 a_{0}\lambda  (n^{2}-2  n+1)}{n  \beta A}} \: ( x \: - \: \frac{A} { 2 \: n \: a_{0}} \: t)))
\end{equation}
so, the travelling wave solution is given by

\begin{equation}
\label{eq:m3st8}
u(x,t)=(\: a_{0}(1+\cot^{2} \: ( \frac{i}{2} \: \sqrt{\frac{a_{0} \: \lambda \: (n^{2}-2n+1)}
       {n \: \beta \: A}} \: ( x-\frac{1}{2} \: \frac{A}{n \: a_{0}} \: t))))^{\frac{-1}{n-1}}
\end{equation}

\bigskip Case (II): consider, $E = \frac{(n^{2} \: a_{2}-2 \: n \: a_{2}+ 2 \: \beta \: n^{2}+ 
                                   2 \:  \beta \: n + a_{2})}{\alpha}$ then,
         $a_{0}=\frac{1}{4} \: \frac{\lambda \: E}{\beta \: n^{2} }$, \; \; $a_{1}=0$, \; $b_{1}=0$, \; $b_{2}=0$, \; $b=\frac{1}{4}\frac{\lambda \: E}{a_{2} \: \beta \: n^{2}}$, \; $c=\frac{a_{2}\: (n^{2}-2 \: n + 1)}{E},$ with $a_{2}$ being an arbitrary constant, then

\begin{equation}
\label{eq:m3st9}
V(x,t) =\frac{1}{4}\frac{E \: \lambda }{\beta \: n^{2}} \: ( 1-\tan ^{2}( 
\frac{i}{2}\sqrt{ \frac{E \: \beta }{\lambda \: n^{2} \: a_{2}}} \:( x-
\frac{a_{2} \: (n^{2}-2n+1)}{E} \: t)))  
\end{equation}
so, the travelling wave solution is given by 
\begin{equation}
\label{eq:m3st10}
u(x,t)=( \frac{1}{4}\frac{E \beta }{\lambda n^{2}} \: ( 1-\tan
       ^{2}( \frac{i}{2}\sqrt{ \frac{E \beta }{\lambda  n^{2} a_{2}} 
       }( x - \frac{a_{2} (n^{2}-2n + 1)}{E} \: t))))^{\dfrac{-1}{n-1}}
\end{equation}

\bigskip Case (III) : $a_{0}=\frac{1}{8} \: \frac{\lambda \: E}{\beta \: n^{2}}$, \;
$a_{1}=0$, \; $b_{1}=0$, \; $b_{2}=\frac{1}{256}\frac{\lambda^{2} \: E}{a_{2} \: \beta^{2}
\: n^{4} \: \alpha}$, \; $b=\frac{1}{16}\frac{\lambda \: E}{\beta  \: n^{2} \: a_{2}}$, \: $c=\frac{a_{2} \: (n^{2}-2n+1)}{E}$ with $a_{2}$ being an arbitrary constant, then
\begin{eqnarray}
\label{eq:m3st11}
  V(x,t) & = & \frac{1}{8}\frac{E \lambda }{\beta n^{2}} \: (1-\frac{1}{2}\tan^{2}( \frac{i}{4} \: \sqrt{\frac{E \: \lambda }{\beta \: n^{2} \: a_{2}}} \: 
( x-\frac{a_{2} \: (n^{2}-2n+1)}{E} \: t)) \nonumber \\
& & {} - \frac{1}{2}\cot^{2}(\frac{i}{4} \: \sqrt{ \frac{E \: \lambda }{\beta \:
n^{2} \: a_{2}}} \: (x-\frac{a_{2} \: (n^{2}-2n+1)}{E} \: t)))
\end{eqnarray}
so, the travelling wave solution is given by

\begin{eqnarray}
\label{eq:m3st12}
  u(x,t) & = & (\frac{1}{16}\frac{ \beta E}{\lambda n^{2}} \: (1-\frac{1}{2}\tan^{2}( \frac{i}{4} \: \sqrt{\frac{E \: \lambda }{\beta \: n^{2} \: a_{2}}} \: 
( x-\frac{a_{2} \: (n^{2}-2n+1)}{E} \: t)) \nonumber \\
& & {} - \frac{1}{2}\cot^{2}(\frac{i}{4} \: \sqrt{ \frac{E \: \lambda }{\beta \:
n^{2} \: a_{2}}} \: (x-\frac{a_{2} \: (n^{2}-2n+1)}{E} \: t)))) ^{\frac{-1}{n-1}}
\end{eqnarray}
\subsection{\protect\bigskip Variant of the PHI-four equation}

\bigskip A second important example is the Variant of the PHI-four equation \cite{soliman:exacttrav}, which can be written as

\begin{equation}
\label{eq:m3nd1}
      u_{tt} - \alpha \: u_{xx} - \lambda \: u + \beta \: u^{n} = 0 \;\;\;  n > 0
\end{equation}

\bigskip where $\alpha$, $\lambda$, and $\beta$ are arbitrary constants. In
order to solve Eq.~\eqref{eq:m3nd1} by the MEDA method, we use the wave transformation 
$u(x,t)= U(z)$ with wave variable $z = i(x - c \: t)$ Eq.~\eqref{eq:m3nd1} takes the 
form of an ordinary differential equation

\begin{equation}
\label{eq:m3nd2}
 c^{2}U^{\prime \prime}-\alpha \: U^{\prime \prime }-\lambda \: U+\beta \: U^{n}= 0
\end{equation}

Or equivalently

\begin{equation}
\label{eq:m3nd3}
      - \lambda \: U + \beta \: U^{n}-(c^{2}-\alpha) \: U^{\prime \prime } = 0
\end{equation}

Balancing the order of $U^{n}$ with the order of $U^{\prime \prime }$ in Eq.
~\eqref{eq:m3nd3}, we find $M = \frac{2}{n-1}$. To obtain a closed form of the analytic solution
we use a transformation formula $U = V^{\frac{1}{n-1}}$, that transforms Eq.~\eqref{eq:m3nd3} to

\begin{equation}
\label{eq:m3nd4}
 -\lambda (n-1)^{2} V^{2}+ \beta (n-1)^{2} V^{3}+(\alpha - c^{2})
   (n-1) V V^{\prime \prime}+( \alpha - c^{2}) (2-n) (V^{\prime })^{2} = 0
\end{equation}

Balancing the order of $V^{3}$ with the order of $V \: V^{\prime \prime}$ in
Eq.~\eqref{eq:m3nd4}, we find $M=2$. So the solution takes the form
\begin{equation}
\label{eq:m3nd5}
      V(z)=a_{0}+ a_{1} \: \phi(z)+ a_{2} \: \phi(z)^{2}+ b_{1} \: \phi
           (z)^{-1}+ b_{2} \: \phi(z)^{-2}
\end{equation}

Substituting Eq.~\eqref{eq:m3nd5} into Eq.~\eqref{eq:m3nd4} and making use of Eq.
~\eqref{eq:m204}, we obtain a system of algebraic equations, for $a_{0}$, $a_{1}$, $a_{2}$, $b_{1}$, $b_{2}$, and $b$. We solve the obtained system of algebraic equations give
the following three cases:

\bigskip Case (I): $a_{0}= \frac{\lambda \: (n+1)}{2 \: \beta }$, $ a_{1}=0$, $ a_{2}=0$, $b_{1}=0$, $b_{2}= \frac{(n+1) \: (n^{2}-2n+1) \: \lambda ^{2}}{8 \: \beta \:
(c^{2}-\alpha )}$, $b=\frac{\lambda \: (n^{2}-2n+1)}{4 \:(c^{2}-\alpha)}$ with $c$
being an arbitrary constant. then

\begin{eqnarray}
\label{eq:m3nd6}
V(x,t)= \frac{\lambda \: (n+1)}{2 \: \beta }+ \frac{\lambda \: (n+1)}{2 \: \beta}
       \: \cot ^{2}\left( \frac{i}{2} \: \sqrt{\frac{\lambda \: (n^{2}-2n+1)}{
       c^{2}-\alpha }} \: (x-c \: t)\right) 
\end{eqnarray}
so, the travelling wave solution is given by
\begin{align} 
\label{eq:m3nd7}
    u(x,t) &= \left( \frac{\lambda (n+1)}{2 \beta }+ \frac{\lambda (n+1)}{
           2 \beta } \: \cot ^{2}\left( \frac{i}{2} \: \sqrt{\frac{\lambda 
           (n^{2}-2n+1)}{c^{2}-\alpha }}(x-c \: t)\right)\right)^{(\frac{1}{n-1})}
\end{align}

\bigskip Case (II): $a_{0}= \frac{\lambda (n+1)}{2 \: \beta }$,\:$a_{1}=0$,
$a_{2}=\frac{2(-\alpha \: n + c^{2} \: n - \alpha + c^{2})}{\beta\:(n^{2}-2n+1)}$,
$b_{1}=0$,\:$b_{2}=0$, $b=\frac{\lambda \: (n^{2}-2\:n+1)}{4\:(c^{2}-\alpha)}$,
with $c$ being an arbitrary constant, then

\begin{eqnarray}
\label{eq:m3nd8}
V(x,t) &=&  \frac{1}{2}\:\frac{\lambda\:(n+1)}{\beta }+ \frac{1}{2}\:\frac{\lambda
\:(-\alpha \:n + c^{2}\: n - \alpha + c^{2})}{\beta \:(c^{2}-\alpha)} \cdot \nonumber \\
& & \tan^{2}\:( \frac{i}{2} \:\sqrt{\frac{\lambda \:(n^{2}-2\:n + 1)}
{c^{2}-\alpha }}\:(x-c\:t))                                        
\end{eqnarray}

so, the travelling wave solution is given by

\begin{eqnarray}
\label{eq:m3nd9}
u(x,t) &=& ( \frac{1}{2}\:\frac{\lambda\:(n+1)}{\beta }+ \frac{1}{2}\:\frac{\lambda
\:(-\alpha \:n + c^{2}\: n - \alpha + c^{2})}{\beta \:(c^{2}-\alpha)} \cdot \nonumber \\
& & \tan^{2}\:( \frac{i}{2} \:\sqrt{\frac{\lambda \:(n^{2}-2\:n + 1)}
{c^{2}-\alpha }}\:(x-c\:t)))^{ \frac{1}{n-1}}
\end{eqnarray}

\bigskip Case(III): Consider, $D = \sqrt{\frac{2\:\lambda \:n -\lambda n^{2}-\lambda
-16 \: \alpha \:b}{b}}$, $a_{0}=\frac{1}{4}\:\frac{\lambda \:(n+1)}{\beta }$, $a_{1}=0$, $a_{2}=\frac{1}{8}\:\frac{(n+1)\: \lambda }{\beta \: b}$, $b_{1}=0$, $b_{2}=\frac{1}{8}\:\frac{(n+1)\:\lambda \:b}{\beta }$, $c=\frac{i}{4}\: D$ 
with $b$ being an arbitrary constant. then,

\begin{eqnarray}
\label{eq:m3nd10}
V(x,t)= \frac{\lambda \: (n+1)}{4\:\beta }+\frac{\lambda \:
(n+1)}{8 \: \beta }\tan^{2}( \sqrt{-b}\:(x-\frac{i}{4}\:D \:t))+  \notag \\
\frac{\lambda \:(n+1)}{8\:\beta }\cot^{2}( \sqrt{-b}\:(x-\frac{i}{4}\:D \:t))
\end{eqnarray}

so, the travelling wave solution is given by
\begin{eqnarray}
\label{eq:m3nd11}
u(x,t)= ( \frac{\lambda \: (n+1)}{4\:\beta }+\frac{\lambda \:
(n+1)}{8 \: \beta }\:\tan^{2}( \sqrt{-b}\:(x-\frac{i}{4}\:D \:t))+  \notag \\
\frac{\lambda \:(n+1)}{8\:\beta }\:\cot^{2}( \sqrt{-b}\:(x-\frac{i}{4}\:D \:t)))^{(\frac{1}{n-1})} 
\end{eqnarray}

\bigskip Case(IV): $a_{0}=\frac{1}{4}\:\frac{\lambda \:(n+1)}{\beta }$, $a_{1}=0$, $a_{2}=\frac{1}{8}\:\frac{(n+1)\: \lambda }{\beta \: b}$, $b_{1}=0$, $b_{2}=\frac{1}{8}\:\frac{(n+1)\:\lambda \:b}{\beta }$, $c= - \frac{i}{4}\: D$ 
with $b$ being an arbitrary constant then,

\begin{eqnarray}
\label{eq:m3nd12}
V(x,t)= \frac{\lambda \: (n+1)}{4\:\beta }+\frac{\lambda \:
(n+1)}{8 \: \beta }\tan^{2}( \sqrt{-b}\:(x+\frac{i}{4}\:D \:t))+  \notag \\
\frac{\lambda \:(n+1)}{8\:\beta }\cot^{2}( \sqrt{-b}\:(x+\frac{i}{4}\:D \:t))
\end{eqnarray}
so, the travelling wave solution is given by
\begin{eqnarray}
\label{eq:m3nd13}
u(x,t)= ( \frac{\lambda \: (n+1)}{4\:\beta }+\frac{\lambda \:
(n+1)}{8 \: \beta }\:\tan^{2}( \sqrt{-b}\:(x+\frac{i}{4}\:D \:t))+  \notag \\
\frac{\lambda \:(n+1)}{8\:\beta }\:\cot^{2}( \sqrt{-b}\:(x+\frac{i}{4}\:D \:t)))^{(\frac{1}{n-1})} 
\end{eqnarray}
All the solutions of the equations are new.
\subsection{The variant Boussinesq equations}

Finally, we consider a very important example as an illustration of the
modified extended direct algebraic method for solving the variant Boussinesq
equations \cite{ali:newexact}, we will consider the following system of equations
\begin{eqnarray}
       u_{t}+v_{x}+u \: u_{x} = 0     \label{eq:m3rd1} \\
       v_{t}+(uv)_{x}+u_{xxx} = 0      \label{eq:m3rd2}
\end{eqnarray}

\qquad\ To solve the system of Eqs~\eqref{eq:m3rd1},~\eqref{eq:m3rd2} by means of the modified
extended direct algebraic method, we use the wave transformation $u(x,t)=U(z)$ and $v(x,t)=V(z)$ with complex wave variable $z = i( x + \lambda \: t)$. Therefore, system \eqref{eq:m3rd1}, \eqref{eq:m3rd2} is reduced to the ordinary differential equations in the form

\begin{eqnarray}
\lambda \: U^{\prime}+V^{\prime }+ \frac{1}{2} \: (U^{\prime })^{2} = 0   \label{eq:m3rd3} \\
\lambda \: V^{\prime}+(UV)^{\prime }-U^{\prime \prime \prime } = 0         \label{eq:m3rd4}
\end{eqnarray}

Integrating both equations once leads to:

\begin{eqnarray}
      C_{1}-\lambda \: U - \frac{1}{2} \: U^{2} = V        \label{eq:m3rd5} \\
      \lambda \: V + U \: V - U^{\prime \prime }  = C_{2}   \label{eq:m3rd6}
\end{eqnarray}

where $C_{1}$ and $C_{2}$ are integrating constants, so as to we find the
special forms of the exact solutions, for simplicity purpose, we take $C_{1}$
$=$$C_{2}=0$ Substituting~\eqref{eq:m3rd5} into~\eqref{eq:m3rd6} gives

\begin{equation}
\label{eq:m3rd7}
    U^{\prime \prime }+\frac{1}{2}U^{3}+\frac{3}{2}\lambda U^{2}+\lambda ^{2}U = 0,
\end{equation}

\bigskip By balancing $U^{\prime \prime }$ with $U^{3}$ in eq.~\eqref{eq:m3rd7}, we find 
$M=1$. So the solutions take the form

\begin{equation}
\label{eq:m3rd8}
      U(z) = a_{0}+a_{1}\phi(z)+ b_{1}\phi (z)^{-1}
\end{equation}

\begin{equation}
\label{eq:m3rd9}
      V(z)=-\lambda \left[ a_{0}+a_{1}\phi(z)+b_{1}\phi(z)^{-1}\right]-\frac{1}{2}
\left[ a_{0}+a_{1}\phi(z)+b_{1}\phi(z)^{-1}\right]^{2}
\end{equation}

Inserting Eqs.~\eqref{eq:m3rd8},~\eqref{eq:m3rd9} into Eq.~\eqref{eq:m3rd6}, making
use of Eq.~\eqref{eq:m204}, and by using Maple Package, we get a system of algebraic
equations, for $a_{0}, a_{1}, b_{1},\lambda$, and $b$. We solve the obtained algebraic 
system of equations by Mable Package and select four cases of solutions as:

\bigskip Case (I): $a_{1}=2\,i$, $b_{1}=0$, $b=\frac{-\:1}{4}\;a_{0}^{2}$, $\lambda
=-\: a_{0}$, with $a_{0}$ being an arbitrary constant, the complex wave solutions are :

\begin{equation}
\label{eq:m3rd10}
      u(x,t)=a_{0}(1+i\tan (\frac{ia}{2}(x-a_{0}t)))
\end{equation}

\begin{equation}
\label{eq:m3rd11}
      v(x,t)=a_{0}(a_{0}(1+i\tan (\frac{ia}{2}(x-a_{0}t))))-
       \frac{1}{2}(a_{0}(1+i\tan (\frac{ia}{2}(x-a_{0}t))))^{2} 
\end{equation}

\bigskip Case (II): $a_{1}= - 2\,i$, $b_{1}=0$, $b=\frac{-\:1}{4}\;a_{0}^{2}$, $\lambda
=-\: a_{0}$, with $a_{0}$ being an arbitrary constant, the complex wave
solutions are :

\begin{equation}
\label{eq:m3rd12}
      u(x,t)=a_{0}(1-i\tan (\frac{ia}{2}(x-a_{0}t)))
\end{equation}

\begin{equation}
\label{eq:m3rd13}
       v(x,t)=a_{0}(a_{0}(1-i\tan (\frac{ia}{2}(x-a_{0}t))))-
       \frac{1}{2}(a_{0}(1-i\tan (\frac{ia}{2}(x-a_{0}t))))^{2} 
\end{equation}

\bigskip Case (III): $a_{1}=2\,i$, $b_{1}=\frac{-\:1}{4}\;a_{0}^{2}\:i$, $b=\frac{-\:1}{8}\:a_{0}^{2}$,
          $\lambda=-\: a_{0}$ with $a_{0}$ being an arbitrary constant, the
complex wave solutions are
\begin{equation}
\label{eq:m3rd14}
u(x,t) = a_{0}+\frac{\sqrt{2}\:a_{0}}{2}(\:\tan\:(\frac{\sqrt{2}\:a_{0}}{4}(x-a_{0}\:t))+
         \:\cot\:(\frac{\sqrt{2}\:a_{0}}{4}(x-a_{0}\:t)))   
\end{equation}
\begin{eqnarray}
\label{eq:m3rd15}
v(x,t)=a_{0}(a_{0}+\frac{\sqrt{2}\:a_{0}}{2}(\tan\:(\frac{\sqrt{2}\:a_{0}}{4}(x-a_{0}\:t))+
       \cot\:(\frac{\sqrt{2}\:a_{0}}{4}(x-a_{0}\:t)))) \notag \\
-\frac{1}{2}(a_{0}(a_{0}+\frac{\sqrt{2}\:a_{0}}{2}(\tan\:(\frac{\sqrt{2}\:a_{0}}{4}
       (x-a_{0}\:t))+ \cot\:(\frac{\sqrt{2}\:a_{0}}{4}(x-a_{0}\:t))))^{2} 
\end{eqnarray}
\bigskip Case (IV):: $a_{1}=-2\,i$, $b_{1}=\frac{1}{4}\;a_{0}^{2}\:i$, $b=\frac{1}{8}\:a_{0}^{2}$,
          $\lambda=-\: a_{0}$ with $a_{0}$ being an arbitrary constant, the
complex wave solution are

\begin{eqnarray}
u(x,t) = a_{0}-\frac{\sqrt{2}\:a_{0}}{2}(\tan\:(\frac{\sqrt{2}\:a_{0}}{4}
         (x-a_{0}\:t))+\cot\:(\frac{\sqrt{2}\:a_{0}}{4}(x-a_{0}\:t)))    \label{eq:m3rd16} \\
v(x,t) = a_{0}( a_{0}-\frac{\sqrt{2}\:a_{0}}{2}(\tan\:(\frac{\sqrt{2}\:a_{0}}{4}
         (x-a_{0}\:t))+\cot\:(\frac{\sqrt{2}\:a_{0}}{4}(x-a_{0}\:t)))) \notag \\
-\frac{1}{2}( a_{0}-\frac{\sqrt{2}\:a_{0}}{2}(\tan\:(\frac{\sqrt{2}\:a_{0}}{4}
         (x-a_{0}\:t))+\cot\:(\frac{\sqrt{2}\:a_{0}}{4}(x-a_{0}\:t))))^{2} \label{eq:m3rd9}
\end{eqnarray}
All the solutions of the variant Boussinesq equations are new.
\section{Conclusions}
In this paper, the MEDA\ method has been successfully applied to find the
solution for three nonlinear partial differential equations such as the the
variant of the RLW\ equation, the variant ~of the PHI-four equation, and the
variant Boussinesq equations . The modified extended direct algebraic method
is used to find a new complex travelling wave solutions. The results show
that the modified extended direct algebraic method is a powerful
mathematical tool to solve the the variant of the RLW\ equation, the variant
~of the PHI-four equation, and the variant Boussinesq equations it is also a
promising method to solve other nonlinear equations.

%

\end{document}